\documentclass[fontsize=12pt,a4paper,headings=normal,
twoside=false,leqno,parskip=half-,abstract=true]{scrartcl}
\usepackage[english]{babel}
\usepackage[utf8]{inputenc}
\usepackage{hyperref}
\hypersetup{
 pdftitle={SturmForQuasilinearEq},
 pdfauthor={Phillipo Lappicy},
 colorlinks=true,
 linkcolor=blue,
 citecolor=blue,
 filecolor=blue,
 urlcolor=blue}

\usepackage{graphicx}
\usepackage[format=plain,labelfont=bf,font=small]{caption}
\usepackage{xcolor}
\usepackage[arrow, matrix, curve]{xy}
\usepackage{float}

\usepackage{caption}
\usepackage{tabulary}
\usepackage{array}
\newcolumntype{N}[1]{>{\centering\arraybackslash}m{#1}}

\usepackage{amsmath,amsthm}
\usepackage{amssymb} 

\makeatletter
\newcommand{\tpitchfork}{%
  \vbox{
    \baselineskip\z@skip
    \lineskip-.52ex
    \lineskiplimit\maxdimen
    \m@th
    \ialign{##\crcr\hidewidth\smash{$-$}\hidewidth\crcr$\pitchfork$\crcr}
  }%
}
\makeatother
\usepackage{latexsym}

\usepackage[notref,notcite,color,final 
]{showkeys}

\definecolor{refkey}{rgb}{1,0,0}
\definecolor{labelkey}{rgb}{1,0,0}

\usepackage{tikz}

\usepackage[textwidth=2cm,textsize=small,backgroundcolor=none]{todonotes}

  \mathchardef\ordinarycolon\mathcode`\:
  \mathcode`\:=\string"8000
  \begingroup \catcode`\:=\active
    \gdef:{\mathrel{\mathop\ordinarycolon}}
  \endgroup

\newtheorem{thm}{Theorem}[section]
\newtheorem{lem}[thm]{Lemma}
\newtheorem{prop}[thm]{Proposition}
\newtheorem{cor}[thm]{Corollary}
\newenvironment{pf}[1][Proof]{\begin{trivlist}
\item[\hskip \labelsep {\bfseries #1}]}{\end{trivlist}}

\usepackage{comment}
\usepackage{caption}

\begin{document}

\title{{\LARGE{Sturm attractors for quasilinear parabolic equations}}}

\author{
 \\
{~}\\
Phillipo Lappicy*\\
\vspace{2cm}}

\date{version of \today}
\maketitle
\thispagestyle{empty}

\vfill

$\ast$\\
Instituto de Ciências Matemáticas e de Computação\\
Universidade de S\~ao Paulo\\
Avenida trabalhador são-carlense 400\\
13566-590, São Carlos, SP, Brazil\\


\newpage
\pagestyle{plain}
\pagenumbering{arabic}
\setcounter{page}{1}

\begin{abstract}
The goal of this paper is to construct explicitly the global attractors of quasilinear parabolic equations, as it was done for the semilinear case by Brunovsk\'y and Fiedler (1986), and generalized by Fiedler and Rocha (1996). In particular, we construct heteroclinic connections between hyperbolic equilibria, stating necessary and sufficient conditions for heteroclinics to occur. Such conditions can be computed through a permutation of the equilibria. Lastly, an example is computed yielding the well known Chafee-Infante attractor.

\ 

\end{abstract}

\section{Main results}\label{sec:intro}

\numberwithin{equation}{section}
\numberwithin{figure}{section}
\numberwithin{table}{section}

Consider the scalar quasilinear parabolic differential equation
\begin{equation}\label{PDEquasi}
    u_t = a(x,u,u_x) u_{xx}+f(x,u,u_x)
\end{equation}
with initial data $u(0,x)=u_0(x)$ such that $f\in C^2$, $a\in C^1$ satisfy the strict parabolicity condition $a(x,u,u_x)\geq\epsilon>0$, and $x \in [0,\pi]$ has Neumann boundary.

The equation \eqref{PDEquasi} defines a semiflow denoted by $(t,u_0)\mapsto u(t)$ in a Banach space $X^\alpha:=C^{2\alpha+\beta}([0,\pi])$. We suppose that $2\alpha+\beta>1$ so that solutions are at least $C^1$. The appropriate functional setting is described in Section \ref{sec:funcquasi}. 

In order to study the long time behaviour of \eqref{PDEquasi}, we suppose that $a$ and $f$ satisfy the following growth conditions, where $p:=u_x$,
\begin{align}\label{dissquasi}
    f(x,u,0) \cdot u&<0\nonumber\\
    |f(x,u,p)|&<f_1(u)+f_2(u)|p|^\gamma\nonumber\\
    \frac{|a_x|}{1+|p|}+|a_u|+|a_p|\cdot [1+|p|]&\leq f_3(|u|)\\
    0<\epsilon \leq a(x,u,p) &\leq \delta \nonumber 
\end{align}
where the first condition holds for $|u|$ large enough, uniformly in $x$, the second for all $(x,u,p)$ for continuous $f_1,f_2$ and $\gamma<2$, the third for continuous $f_3$ and $\epsilon,\delta>0$. 

Those conditions imply that $|u|$ and $|u_x|$ are bounded. Hence, the semiflow is dissipative: trajectories $u(t)$ eventually enter a large ball in the phase-space $X^\alpha$. See Chapter 6, Section 5 in \cite{Ladyzhenskaya68}. Also \cite{KruzhkovOleinek61} and \cite{BabinVishik92}. 

Moreover, these hypothesis guarantee the existence of a nonempty global attractor $\mathcal{A}$ of \eqref{PDEquasi}, which is the maximal compact invariant set. Equivalently, it is the set of bounded trajectories $u(t)$ in the phase-space $X^\alpha$ that exist for all $t\in \mathbb{R}$. See \cite{BabinVishik92}. A new perspective of solutions which do not exist for all times is given by Stuke \cite{Stuke17}.

The goal of this paper is to decompose $\mathcal{A}$ into smaller invariant sets, and describe how those sets are related.

Due to a Lyapunov functional, constructed by Matano \cite{Matano88} and Zelenyak \cite{Zelenyak68}, the global attractor consists of equilibria and their heteroclinic connections within their unstable manifolds. It still persits the question of which equilibria connects to which other. This geometric description was carried out by Hale and do Nascimento \cite{HaleNascimento83} for the Chafee Infante problem, by Brunovský and Fiedler \cite{FiedlerBrunovsky89} for $f(x,u)$ and by Fiedler and Rocha \cite{FiedlerRocha96} for $f(x,u,u_x)$. Such attractors are known as \emph{Sturm attractors}.

Constructing the Sturm attractor for the equation \eqref{PDEquasi} is problematic due to its quasilinear nature. It is the aim of this chapter to modify the existing theory for such equations and still obtain a Sturm attractor. 

For the statement of the main theorem of this chapter, denote by the \emph{zero number} $z(u_*)$ the number of sign changes of a continuous function $u_*(x)$. Recall that the \emph{Morse index} $i(u_*)$ of an equilibrium $u_*$ is given by the number of positive eigenvalues of the linearized operator at such equilibrium, that is, the dimension of the unstable manifold of said equilibrium.

Recall an equilibrium $u_*(x)$ is \emph{hyperbolic} if the linearization operator of the right hand of \eqref{PDEquasi} has no eigenvalue in the imaginary axis.

We say that two different equilibria $u_-,u_+$ of \eqref{PDEquasi} are \emph{adjacent} if there does not exist an equilibrium $u_*$ of \eqref{PDEquasi} such that $u_*(0)$ lies between $u_-(0)$ and $u_+(0)$, and
\begin{equation*}
    z(u_--u_*) = z(u_--u_+) = z(u_+-u_*).
\end{equation*}
This notion was firstly described by Wolfrum \cite{Wolfrum02}.

Both the zero number and Morse index can be computed from a permutation of the equilibria, as it was done in \cite{FuscoRocha} and \cite{FiedlerRocha96}. Such permutation is called the \emph{Sturm Permutation} and is computed in Section \eqref{sec:permquasi}, as it was done \cite{FiedlerRocha96}. For such, it is required that the flow of the equilibria equation of \eqref{PDEquasi} exists for all $x\in [0,\pi]$. 

\begin{thm}\emph{\textbf{Sturm Attractor} } \label{attractorthmquasi}

Consider $a\in C^1$ and $f\in C^2$ satisfying the growth conditions \eqref{dissquasi}. Suppose that all equilibria for the equation \eqref{PDEquasi} are hyperbolic. Then, 
\begin{enumerate}
    \item the global attractor $\mathcal{A}$ of \eqref{PDEquasi} consists of finitely many equilibria $\mathcal{E}$ and heteroclinic connections $\mathcal{H}$.
    \item there is a heteroclinic $u(t)\in \mathcal{H}$ between two equilibria $u_-,u_+\in\mathcal{E}$ so that
    \begin{equation*}
        u(t)\to_{t\to \pm\infty} u_{\pm}
    \end{equation*}
    if, and only if, $u_-$ and $u_+$ are adjacent and $i(u_-)>i(u_+)$.
\end{enumerate}
\end{thm}

This meets some expectations of Fiedler \cite{Fiedler96}, which mentions that fully nonlinear equations yield the same type of attractors as the semilinear ones. We prove it for the quasilinear, and leave the fully nonlinear for another occasion.

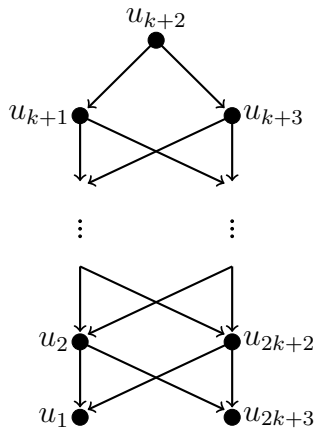
\begin{figure}[ht]\centering
\begin{tikzpicture}
\filldraw [black] (0,0) circle (3pt) node[anchor=south]{$u_{k+2}$};
\filldraw [black] (-1,-1) circle (3pt) node[anchor=east]{$u_{k+1}$};
\filldraw [black] (1,-1) circle (3pt) node[anchor=west]{$u_{k+3}$};

\filldraw [black] (-1,-2.4) circle (0.5pt);
\filldraw [black] (-1,-2.5) circle (0.5pt);
\filldraw [black] (-1,-2.6) circle (0.5pt);
\filldraw [black] (1,-2.4) circle (0.5pt);
\filldraw [black] (1,-2.5) circle (0.5pt);
\filldraw [black] (1,-2.6) circle (0.5pt);

\filldraw [black] (-1,-4) circle (3pt) node[anchor=east]{$u_{2}$};
\filldraw [black] (1,-4) circle (3pt) node[anchor=west]{$u_{2k+2}$};
\filldraw [black] (-1,-5) circle (3pt) node[anchor=east]{$u_{1}$};
\filldraw [black] (1,-5) circle (3pt) node[anchor=west]{$u_{2k+3}$};

\draw[thick,->] (0,0) -- (-0.9,-0.9);
\draw[thick,->] (0,0) -- (0.9,-0.9);

\draw[thick,->] (-1,-1) -- (0.9,-1.9);
\draw[thick,->] (1,-1) -- (-0.9,-1.9);
\draw[thick,->] (-1,-1) -- (-1,-1.87);
\draw[thick,->] (1,-1) -- (1,-1.87);

\draw[thick,->] (-1,-3) -- (0.9,-3.9);
\draw[thick,->] (1,-3) -- (-0.9,-3.9);
\draw[thick,->] (-1,-3) -- (-1,-3.87);
\draw[thick,->] (1,-3) -- (1,-3.87);

\draw[thick,->] (-1,-4) -- (0.9,-4.9);
\draw[thick,->] (1,-4) -- (-0.9,-4.9);
\draw[thick,->] (-1,-4) -- (-1,-4.87);
\draw[thick,->] (1,-4) -- (1,-4.87);

\end{tikzpicture}
\caption{Global attractor $\mathcal{A}$ of Chafee-Infante type} \label{FIGCOR}
\end{figure}

In particular, we compute and comment on some explicit attractors. Firstly, when the diffusion and reaction are balanced with adjusted powers, and of the Chafee Infante type. This attractor could be used as an application of the Einstein Hamiltonian equation, when there is still a degenerate term on the diffusion operator. Such problem shall be treated soon. See \cite{LappicyPhD}. This same diffusion with another reaction is used to model the curve shortening flow in $\mathbb{R}^2$, as in \cite{Angenent91}. See the discussion section below for further exploration on such topics.
\begin{cor}\emph{\textbf{Chafee-Infante Attractor} } \label{CIquasi}

Consider the equation \eqref{PDEquasi} with $f(\lambda, u)=\lambda a(x,u,u_x)u(1-u^2)$. Let $\lambda\in(\lambda_k,\lambda_{k+1})$, where $\lambda_k$ is the $k$-th eigenvalue of the Laplacian with $k\in\mathbb{N}_0$. 

Then, there are $2k+3$ hyperbolic equilibria $u_1,...,u_{2k+3}$ and the connections between equilibria within the attractor $\mathcal{A}$ is described in the Figure \ref{FIGCOR}.
\end{cor}

The remaining is organized as follows. 

We firstly introduce the functional setting in Section \ref{sec:funcquasi}, including invariant manifolds and Matano's Lyapunov functional for the quasilinear case. In particular this implies that the attractor consists of equilibria and heterolinics. we also recall the dropping lemma and some consequences, which hold since the difference of solutions satisfy a linear equation. This is a fundamental result for the attractor construction that dates back to Sturm \cite{Sturm}.

Then, we focus on the connection problem. All the necessary information about the adjacency, namely the zero numbers and Morse indices, are encoded in a permutation of the equilibria, which is described in Section \ref{sec:permquasi}. This was done firstly by \cite{FuscoRocha}. The shooting is similar to the semilinear case, but one needs to divide the reaction term by the diffusion coefficient $a$, noticing $a>0$. All the previous tools are put together to construct the attractor in Section \ref{sec:globalsturmquasi}, as it was done \cite{FiedlerRocha96}. 

Lastly, Section \ref{sec:CIquasi} gives an example of the developed theory and constructs the attractor of the Chafee-Infante type, and a discussion is carried in \ref{sec:discuss}, regarding applications and generalizations of the result.

\section{Proof of main result}

\subsection{Background} \label{sec:funcquasi}





The Banach space used on the upcoming theory consists on subspaces of H\"older continuous functions $C^\beta([0,\pi])$ with $\beta\in (0,1)$. A more precise description is given below, following \cite{Lunardi95}, \cite{Amann95}, \cite{BabinVishik92}. The notation $C^{\beta}$ for some $\beta\in\mathbb{R}$ indicates that $\beta=[\beta]+\{\beta\}$, where the integer part $[\beta]\in\mathbb{N}$ denotes the $[\beta]$-times differentiable functions whose $[\beta]$-derivatives are $\{\beta\}$-Hölder, where $\{\beta\}\in[0,1)$ is the fractional part of $\beta$.

The equation \eqref{PDEquasi} can be seen as an abstract differential equation on a Banach space,
\begin{equation}\label{linequiv}
    u_t=Au+g(u)    
\end{equation}
where $A:D(A)\rightarrow X$ is the linearization of the right-hand side of \eqref{PDEquasi} at the initial data $u_0(x)$, and the Nemitskii operator $g$ of the remaining nonlinear part, which takes values in $X$, namely $g(u):=a(x,u,u_x)u_{xx}+f(x,u,u_x)-Au$. The spaces considered are $X:=C^\beta([0,\pi])$, and $D(A)=C^{2+\beta}([0,\pi])\subset X$ is the domain of the operator $A$, where $\beta\in (0,1)$. For such setting, see the Chapter 8, in particular the proof of Theorem 8.1.1 in Lunardi \cite{Lunardi95}.


As in Lunardi \cite{Lunardi95}, we consider the interpolation spaces $X^\alpha=C^{2\alpha+\beta}([0,\pi])$ between $D(A)$ and $X$ with $\alpha\in (0,1)$ such that $A$ generates a strongly continuous semigroup in $X^\alpha$, and hence the equation \eqref{PDEquasi} with the dissipative conditions \eqref{dissquasi} defines a dissipative dynamical system in $X^\alpha$. We suppose that $2\alpha+\beta>1$ so that solutions are at least in $C^1([0,\pi])$.

In particular, it settles the theory of existence and uniqueness. For certain qualitative properties of solutions, such as the existence of invariant manifolds tangent to the linear eigenspaces, one needs to know the spectrum of $A$.


Note that $Au=\lambda u$ is a regular Sturm-Liouville problem, since the coefficients depend only on $x$ and are all bounded. Therefore, the spectrum $\sigma (A)$ consists of real simple eigenvalues $\lambda_k$ accumulating at $-\infty$, and corresponding eigenfunctions $\phi_k(x)$ which form a an orthonormal basis of $X$. Moreover, there is a spectral gap between eigenvalues that allows us to get the following filtration of invariant manifolds. Note that since it is supposed that equilibria are hyperbolic, then there is no 0 as an eigenvalue and no center direction.

\begin{prop} \label{hierarchyquasi}
    \emph{\textbf{Filtration of Invariant Manifolds} \cite{Mielke91}, \cite{Lunardi95}}
    
    Let $u_*$ be a hyperbolic equilibrium of \eqref{PDEquasi} with Morse index $n:=i(u_*)$. Then there exists a filtration of the \emph{unstable manifold}
    \begin{equation*}
        W^u_0(u_*)\subset ... \subset W^u_{n-1}(u_*)=W^u(u_*)
    \end{equation*}
    where each $W^u_k$ has dimension $k+1$ and tangent space at $u_*$ spanned by $\phi_0,...,\phi_{k}$. 
    
    Analogously, there is a filtration of the \emph{stable manifold}
    \begin{equation*}
        ... \subset W^s_{n+1}(u_*)\subset W^s_{n}(u_*)=W^s(u_*)
    \end{equation*}
    where each $W^s_k$ has codimension $k$ and tangent space spanned by $\phi_{k},\phi_{k+1},...$.
\end{prop}

Note that the above index labels are not in agreement with the dimension of each submanifold within the filtration, but it is with the number of zeros the corresponding eigenfunction has. For example, each eigenfunction $\phi_k$ has $k$ simple zeroes, whereas the $\dim W^u_k=k+1$. 

An important property is the behaviour of solutions within each submanifold of the above filtration of the unstable or stable manifolds.
\begin{prop} \label{ConvEFquasi}
    \emph{\textbf{Linear Asymptotic Behaviour} \cite{Henry81},  \cite{Angenent86}, \cite{FiedlerBrunovsky86}} 
    
    Consider a hyperbolic equilibrium $u_*$ with Morse index $n:=i(u_*)$ and a trajectory $u(t)$ of \eqref{PDEquasi}. The following holds,
    \begin{enumerate}
    \item If $u(t)\in W^u_k(u_*)\backslash W^u_{k-1}(u_*)$ with $k=0,...,i(u_*)-1$. Then,
        \begin{equation*}
        \frac{u(t)-u_*}{||u(t)-u_*||}\xrightarrow{t\rightarrow -\infty} \pm \phi_k
        \end{equation*}
    \item If $u(t)$ in $W^s_k(u_*)\backslash W^s_{k+1}(u_*)$ with $k\geq i(u_*)$. Then,
        \begin{equation*}
        \frac{u(t)-u_*}{||u(t)-u_*||}\xrightarrow{t\rightarrow \infty} \pm \phi_k
        \end{equation*}
    \end{enumerate}

    where the convergence takes place in $C^1$. 
    
    The conclusions of 1. and 2. also hold true by replacing the difference $u(t)-u_*$ with the tangent vector $u_t$.
\end{prop}

The reason this theorem works for both the tangent vector $v:=u_t$ or the difference $v:=u_1-u_2$ of any two solutions $u_1,u_2$ of the nonlinear equation \eqref{PDEquasi} is that they satisfy a linear equation of the type
\begin{equation}\label{linPDEquasi}
    u_t=a(t,x)u_{xx}+b(t,x)u_x+c(t,x)u
\end{equation}
where $x\in(0,\pi)$ has Neumann boundary conditions and the functions $a(t,x),b(t,x)$ and $c(t,x)$ are bounded. 

The proof in \cite{Angenent86} works for the case $a=a(x,u,u_x)$ and $f=(x,u,u_x)$ considering Dirichlet boundary condition; 
or in case $a=a(x,u)$ and $f=f(x,u)$ considering other boundary conditions. 
For the general case, see \cite{FiedlerBrunovsky86}.







Lastly, we recall the existence of a Lyapunov functional, as it was done by Matano \cite{Matano88}. Hence, bounded trajectories tend to equilibria.

\begin{lem}
    \emph{\textbf{Lyapunov Functional} \cite{Matano88}}
    
    There exists a Lagrange functional $L$ such that 
    \begin{equation}\label{Lyapquasi}
        E:= \int_{0}^\pi L(x,u,u_x) dx
    \end{equation}
    is a Lyapunov functional for the equation \eqref{PDEquasi}.
\end{lem}

Matano's idea yields a Lyapunov functional of the type
\begin{equation}\label{BOOK}
    \frac{dE}{dt}:= -\int_{0}^\pi \frac{L_{pp}(x,u,p)}{a(x,u,p)} (u_t)^2 dx\leq 0
\end{equation}
where $p:=u_x$ and $L$ satisfy the convexity condition $L_{pp}>0$. Hence, one needs that $a(x,u,p)>0$.

Therefore, the LaSalle invariance principle holds and implies that bounded solutions converge to equilibria, and any $\omega$-limit set consists of a single equilibrium. See \cite{Matano88}. Moreover, the global attractor can be characterized as follows, yielding the first part of the main result.

\begin{prop}\emph{\textbf{Attractor Decomposition} \cite{BabinVishik92}}

If all equilibria of \eqref{PDEquasi} are hyperbolic, then the global attractor $\mathcal{A}$ is decomposed as
    \begin{equation*}
        \mathcal{A}=\bigcup_{v\in \mathcal{E}} W^u(v)
    \end{equation*}
    and consists only of the set of finitely many equilibria $\mathcal{E}$ and connection orbits $\mathcal{H}$ within the unstable manifolds.
\end{prop}

Note that hyperbolic equilibria must be isolated. Moreover, there must be finitely many due to dissipativity. This settles the first part of the Theorem \ref{attractorthmquasi}, and we now continue the proof of the second part of such theorem, yielding heteroclinic connections.

Lastly, we present an important tool for the proofs and some consequences.

Let the \emph{zero number} $0\leq z^t(u)\leq \infty$ count the number of strict sign changes in $x$ of a $C^1$ function $u(t,x)\not \equiv 0$, for each fixed $t$. More precisely,
\begin{equation*}
    z^t(u):= \sup_k \left\{  
        \begin{array}{c} 
        \text{There is a partition $\{ x_j\}_{j=1}^{k}$ of } [0,\pi]\\
        \text{such that } u(t,x_j)u(t,x_{j+1})<0 \text{ for all } j
        \end{array} \right\}
\end{equation*}
and $z^t(u)=-1$ if $u\equiv 0$. In case $u$ does not depend on $t$, we omit the index and simply write $z^t(u)=z(u)$.

A point $(t_0,x_0)\in\mathbb{R}\times [0,\pi]$ such that $u(t_0,x_0)=0$ is said to be a \emph{simple zero} if $u_x(t_0,x_0)\neq 0$ and a \emph{multiple zero} if $u_x(t_0,x_0)=0$.

The following result proves that the zero number of certain solutions of \eqref{PDEquasi} is nonincreasing in time $t$, and decreases whenever a multiple zero occur. Different versions of this well known fact are due to Sturm \cite{Sturm}, Matano \cite{Matano82}, Angenent \cite{Angenent88} and others. 
\begin{lem} \label{droplemquasi}
    \emph{\textbf{Dropping Lemma}}
    
    Consider $u\not \equiv 0$ a solution of the linear equation \eqref{linPDEquasi} for $t\in [0,T)$. Then, its zero number $z^t(u)$ satisfies
    \begin{enumerate}
    \item $z^t(u)<\infty$ for any $t\in (0,T)$.
    \item $z^t(u)$ is nonincreasing in time $t$.
    \item $z^t(u)$ decreases at multiple zeros $(t_0,x_0)$ of $u$, that is, 
    \begin{equation*}
        z^{t_0-\epsilon}(u)>z^{t_0+\epsilon}(u)
    \end{equation*} 
    for any  sufficiently small $\epsilon>0$.
    \end{enumerate}
\end{lem}

Recall that both the tangent vector $u_t$ and the difference $u_1-u_2$ of two solutions $u_1,u_2$ of the nonlinear equation \eqref{PDEquasi} satisfy a linear equation as \eqref{linPDEquasi}, and the proof is exactly the same as the one for the semilinear equations. The proof even holds for fully nonlinear equations. 

We mention two consequences of the dropping lemma \ref{droplemquasi} and the asymptotic description in Proposition \ref{ConvEFquasi}. The first is a result relating the zero number within invariant manifolds and the Morse indices of equilibria. The second is the Morse-Smale property.
\begin{prop} \label{Znuminvmfldquasi}
    \emph{\textbf{Zero number and Invariant Manifolds} \cite{FiedlerBrunovsky86} }
    
    Consider an equilibrium $u_*\in\mathcal{E}$ and a trajectory $u(t)$ of \eqref{PDEquasi}. Then,
    \begin{enumerate}
    \item If $u(t) \in W^u(u_*)$, then $i(u_*)>z^t(u-u_*)$.
    \item If $u(t) \in W^s_{loc}(u_*)\backslash \{u_*\}$, then $z^t(u-u_*)\geq i(u_*)$.
    \end{enumerate}
    These results also hold by replacing $u(t)-u_*$ with the tangent vector $u_t$.
\end{prop}

The above theorem implies that \eqref{PDEquasi} has no homoclinic orbits. Indeed, if there were any, then $i(u_*)<i(u_*)$, which is a contradiction.
    

This last theorem implies that if the semigroup has a finite number of equilibria, in which all are hyperbolic, then it is a Morse-Smale system in the sense of \cite{HaleMagalhaesOliva84}. 

\subsection{Sturm permutation}\label{sec:permquasi}

The next step on our quest to find the Sturm attractor is to construct a permutation associated to the equilibria, which is done using shooting methods. This enables the computation of the Morse indices and zero number of equilibria. That was firstly done by Fusco and Rocha \cite{FuscoRocha} using methods also described by Fusco, Hale and Rocha in \cite{Rocha85}, \cite{RochaHale85}, \cite{Rocha88}, \cite{Rocha94} and \cite{FuscoHale85}.

The equilibria equation associated to \eqref{PDEquasi} can be rewritten as
\begin{equation*}
    0=a(x,u,u_x)u_{xx}+f(x,u,u_x) 
\end{equation*}
for $x \in [0,\pi]$ with Neumann boundary conditions and the parabolicity condition $a>0$.

Reduce the equation to a first order system through $p:=u_\tau$, adding the extra equation $x_\tau=1$ to obtain an autonomous system. Hence, 
\begin{align}\label{shootflowquasi}
\begin{cases}
    u_{\tau}&= p\\
	p_{\tau}&=- \frac{f(x,u,p)}{a(x,u,p)}\\
	x_\tau&=1
\end{cases}
\end{align}
where the Neumann boundary condition becomes 
$p=0$. We suppose that solutions are defined for all $x\in[0,\pi]$ and any initial data.

The idea to find equilibria \eqref{PDEquasi} is as follows. They must lie in the line 
\begin{equation*}
    L_0:=\{(x,u,p)\in\mathbb{R}^3 \textbf{ $|$ }  (x,u,p)=(0,b,0) \text{ and } b\in\mathbb{R}\}  
\end{equation*}
due to Neumann boundary at $x=0$. Then, evolve this line under the flow of the equilibria differential equation and intersect it with an analogous line $L_\pi$ at $x=\pi$, so that it also satisfies Neumann at $x=\pi$. More precisely, one can write the \emph{shooting manifold} as
\begin{equation*}
    M:= \{(x,u,p)\in[0,\pi]\times \mathbb{R}^3 \textbf{ $|$ } (x,u(x,b,0),p(x,b,0)),b\in\mathbb{R}  \}.
\end{equation*}
where $(x,u(x,b,0),p(x,b,0))$ is the solution of \eqref{shootflowquasi} which evolves the initial data $(0,b,0)$.

Denote by $M_x$ the cross-section of $M$ for some fixed $x\in [0,\pi]$. This is a curve parametrized by $b\in\mathbb{R}$.

We obtain the following characterization of equilibria through the shooting manifolds and its relation with the Morse indices and zero numbers, similar to \cite{Rocha85} and \cite{Hale99}.
\begin{lem}
    \textbf{Equilibria Through Shooting}
    \begin{enumerate}
    \item The set of equilibria $\mathcal{E}$ of \eqref{PDEquasi} is in one-to-one correspondence with $M_{\pi}\cap L_{\pi}$.
    \item An equilibrium point corresponding to fixed $b\in\mathbb{R}$ is hyperbolic if, and only if, $M_\pi$ intersects $L_\pi$ transversely at $(\pi,u(\pi,b,0),0)$.
    \item  If $u_*$ correspond to a hyperbolic equilibrium of \eqref{PDEquasi}, then its Morse index is given by $i(u_*)=1+\lfloor\frac{\zeta(x_0)}{\pi}\rfloor$ where $\zeta(x_0)$ is the angle between $M_\pi$ and $L_\pi$ measured clockwise at their intersection point $x_0$, and $\lfloor.\rfloor$ denotes the floor function. 
    \end{enumerate}
\end{lem}
\begin{pf}
To prove 1), note that a point in $M_{\pi}\cap L_{\pi}$ satisfies the equilibria equation with Neumann boundary conditions by definition of the shooting manifolds. Conversely, consider an equilibrium of \eqref{PDEquasi} satisfying Neumann boundary must be in $M_{\pi}\cap L_{\pi}$. Due to the uniqueness of the shooting differential equation \eqref{shootflowquasi}, such correspondence above is one-to-one. 

To prove 2), consider an equilibrium $u_*$ corresponding to the initial data $b\in\mathbb{R}$. We compare the eigenvalue problem for $u_*$ and the differential equation satisfied by the angle of the tangent vectors of the shooting manifold. 

The eigenvalue problem for $u_*$ is
\begin{equation*}
    \lambda u =a_*(x)u_{xx}+b_*(x)u+ c_*(x) u_x
\end{equation*}
where $x\in[0,\pi]$ has Neumann boundary conditions, and the coefficients are
\begin{align*}
    a_*(x)&:=a(x,u_*,p_*)\\ b_*(x)&:=a_u(x,u_*,p_*).(u_*)_{xx}+D_uf(x,u_*,p_*)\\ c_*(x)&:=a_p(x,u_*,p_*).(u_*)_{xx}+D_pf(x,u_*,p_*).
\end{align*}

Rewriting the above equation as a first order system through $p:=u_x$,
\begin{align*}
\begin{cases}
    u_x&=p\\
    p_x&=-\frac{b_*u+c_*p-\lambda u}{a_*}
\end{cases}
\end{align*}
with Neumann boundary conditions.

In polar coordinates $(u,p)=:(r\cos(\mu),-r\sin(\mu))$, the angle given by $\mu:=\arctan(\frac{p}{u})$ satisfies 
\begin{equation}\label{EVpolarquasi}
    \mu_\tau= \sin^2(\mu)+ \frac{b_*u+c_*p-\lambda u}{a_*} \cos^2(\mu)
\end{equation}
with $\mu(0)=0$ and $\mu(\pi)=k\pi$ for some $k\geq0$. 

On the other hand, $M_x$ is parametrized by the initial data $b\in\mathbb{R}$ and its tangent vector $(\frac{\partial u(x,b)}{\partial b},\frac{\partial p(x,b)}{\partial b})$ corresponding to the trajectory $u_*$ satisfies the following linearized equation, 
\begin{align}\label{tangentshootquasi}
\begin{cases}
    (u_{b})_{x}&= p_{b}\\
	(p_{b})_{x}&=-\frac{b_* u_{b}+ c_* p_{b}}{a_*}
\end{cases}
\end{align}
with initial data $ (u_{b}(0),p_{b}(0))=(1,0)$.

In polar coordinates $(u_{b},p_{b})=:({\rho}\cos({\nu}),-{\rho}\sin({\nu}))$, where ${\nu}$ is the clockwise angle of the tangent vector of $M_x$ at the trajectory $u_*$ with the $u$-axis, 
\begin{equation}\label{tangentpolarquasi}
    {\nu}_x= \sin^2({\nu})+ \frac{b_* u_{b}+ c_* p_{b}}{a_*} \cos^2({\nu})
\end{equation}
with initial data $ {\nu}(0,b,0)=0$. 

Note that the angle $\nu$ of the tangent vector in \eqref{tangentpolarquasi} satisfy is the same equation as the eigenvalue problem in polar coordinates \eqref{EVpolarquasi} with $\lambda=0$ and same boundary conditions at $x=0$. 

Suppose that $u_*$ is not hyperbolic, that is, there exists a solution of \eqref{EVpolarquasi} with $\mu(\pi)=k\pi$ for $\lambda=0$ and some $k\in\mathbb{N}$. Since this is the same equation as \eqref{tangentpolarquasi}, uniqueness implies that $\nu(\pi)=k\pi$. This implies that $M_\pi$ and $L_\pi$ are not transverse.

Conversely, if $M_\pi$ and $L_\pi$ are not transverse, then $\nu(\pi)=k\pi$ for some $k\in\mathbb{N}$. Again, notice this is the same equation for \eqref{EVpolarquasi} and hence there exists a solution of \eqref{EVpolarquasi} for $\lambda=0$ such that $\nu(\pi)=k\pi$. Hence, $\lambda=0$ is an eigenvalue and $u_*$ is not hyperbolic.

To prove 3), 
consider the solution $\mu(x,\lambda)$ of the eigenvalue problem in polar coordinates \eqref{EVpolarquasi}. The Sturm oscillation theorem implies that 
\begin{equation*}
    \psi(\lambda):=\mu(\pi,\lambda)
\end{equation*} 
is decreasing so that $\lim_{\lambda\rightarrow -\infty} \psi(\lambda)=\infty$ and $\lim_{\lambda\rightarrow \infty} \psi(\lambda)=-\pi/2$. Hence, there exists a decreasing sequence $\{ \lambda_k \}_{k\in N}$ to $-\infty$ such that $\psi(\lambda_k)=k\pi$ for $k\in \mathbb{N}$. This implies that there exists a solution of \eqref{EVpolarquasi} for each $\lambda_k$ such that $\psi(\lambda_k)=k\pi$, and hence $\{ \lambda_k \}_{k\in N}$ are the eigenvalues.

Recall that the Morse index $i(u_*)$ is the number of positive eigenvalues of the linearization at $u_*$, that is 
\begin{equation*}
    ...<\lambda_{i(u_*)}<0<\lambda_{i(u_*)-1}<...<\lambda_0. 
\end{equation*}

Since $\psi(\lambda)$ is decreasing and $\lambda_{i(u_*)}$ are eigenvalues, then
\begin{equation*}
    i(u_*)\pi=\psi(\lambda_{i(u_*)})> \psi(0)> \psi(\lambda_{i(u_*)-1})=(i(u_*)-1)\pi .
\end{equation*}

Divide the above by $\pi$ and consider the integer value, yielding that $i(u_*)=\lfloor \frac{\psi(0)}{\pi}\rfloor +1$. By definition, $\psi(0)=\nu(\pi,0)$, which is exactly the angle between $M_\pi$ and $L_\pi$.
\begin{flushright}
	$\blacksquare$
\end{flushright}
\end{pf}

Therefore, a \emph{Sturm permutation} $\sigma$ is obtained by labeling the intersection points $u_i\in M_{\pi}\cap L_\pi$ firstly along $M_{\pi}$ following its parametrization given by $(\pi,u(\pi,b,0),p(\pi,b,0))$ as $b$ goes from $-\infty$ to $\infty$,
\begin{equation*}
    u_1 < ... < u_N
\end{equation*} 
where $N$ denotes the number of equilibria. Secondly, label the intersection points along $L_\pi$ by increasing values,
\begin{equation*}
    u_{\sigma(1)} < ... < u_{\sigma(N)}
\end{equation*} 

The Morse indices of equilibria and the zero number of difference of equilibria can be calculated through the Sturm permutation $\sigma$, as in \cite{Rocha91} and \cite{FiedlerRocha96}. The main tool for such proofs is the third part of the above Lemma: the positive rotation along the shooting curve increases the Morse index.

\subsection{Sturm global structure}\label{sec:globalsturmquasi}

This section gathers all the tools developed in the previous sections, in order to construct the attractor for the quasilinear parabolic equation \eqref{PDEquasi} and prove the second part of the main Theorem \ref{attractorthmquasi}.

Its proof is a consequence of the following three propositions. Firstly, due to the \emph{cascading principle}, it is enough to construct all heteroclinics between equilibria such that the Morse indices of such equilibria differ by 1. Secondly, on one direction, the \emph{blocking principle}: some conditions imply that there does not exist a heteroclinic connection; on the other direction, the \emph{liberalism principle}: if those conditions are violated, then there exists a heteroclinic. Thirdly, \emph{Wolfrum's result} yield a relation between two notion of adjacencies: one that depends on a cascade between equilibria, and one that does not depend on a cascade.

The cascading (Lemma 1.5 in \cite{FiedlerRocha96}) and blocking principles follow from the dropping lemma and its consequences, as Fiedler and Rocha \cite{FiedlerRocha96}. There is only a mild modification in the proof of the liberalism principle (Lemma 1.7 in \cite{FiedlerRocha96}).

\begin{prop}\emph{\textbf{Cascading Principle} \cite{FiedlerRocha96}}\label{cascadingquasi}
 
 There exists a heteroclinic between two equilibria $u_\pm$ such that  $n:=i(u_-)-i(u_+)>0$ if, and only if, there exists a sequence (cascade) of equilibria $\{ v_k\}_{k=0}^n$ with $v_0:=u_-$ and $v_n:=u_+$, such that the following holds for all $k=0,...,n-1$
 \begin{enumerate}
     \item $i(v_{k+1})=i(v_k)+1$
     \item There exists a heteroclinic from $v_{k+1}$ to $v_k$
 \end{enumerate} 
\end{prop}

\begin{prop} \emph{\textbf{Blocking and Liberalism Principles} \cite{FiedlerRocha96}} \label{estabhetsquasi}

There exists a heteroclinic between the equilibria $v_{k+1}$ and $v_k$ with $i(v_{k+1})=i(v_k)+1$ if, and only if,
  \begin{enumerate}
     \item \emph{Morse permit: } $z(v_{k+1}-v_k)=i(v_k)$,
     \item \emph{Zero number permit: } it does not exist an equilibrium $u_*$ between $v_{k+1}$ and $v_k$ at $x=0$ such that $z(v_{k+1}-u_*)=z(v_k-u_*)$.
 \end{enumerate} 
\end{prop}

The blocking and liberalism principles assert that the Morse indices $i(.)$ and zero numbers $z(.)$ construct the global structure of the attractor explicitly. Those numbers can be obtained from the Sturm permutation, as in Section \ref{sec:permquasi}. 

The two propositions above yield the existence of heteroclinics between $u_-$ and $u_+$ if they are \emph{cascadly adjacent}, namely, if there exists a cascade of equilibria $\{ v_k\}_{k=0}^n$ with $v_0:=u_-$ and $v_n:=u_+$ such that for all $k=0,...,n-1$ the following three conditions hold:
  \begin{enumerate}
     \item $i(v_{k+1})=i(v_k)+1$,
     \item $z(v_k-v_{k+1})=i(v_{k+1})$,
     \item there does not exists an equilibrium $v_*$ between $v_k$ and $v_{k+1}$ at $x=0$ satisfying $z(v_k-v_*)=z(v_{k+1}-v_*)$.
 \end{enumerate} 
 
On the other hand, the main Theorem \ref{attractorthmquasi} yields a result through the notion of adjacency in the introduction, which only does not involve a cascade. These notions of adjacency coincide, and this is the core of Wolfrum's ideas in \cite{Wolfrum02}. 

\begin{prop} \emph{\textbf{Wolfrum's equivalence}} \label{wolfrumlemma}
Consider two equilibria $u_\pm\in\mathcal{E}$ such that $n:=i(u_-)-i(u_+)>0$. The equilibria $u_\pm$ are adjacent if, and only if they are cascadly adjacent.
\end{prop}
For the proof of the liberalism theorem, it is used the Conley index to detect orbits between $u_-$ and $u_+$. We give a brief introduction of Conley's theory, and how it can be applied in this context. See Chapters 22 to 24 in \cite{Smoller83} for a brief account of the Conley index, and its extension to infinite dimensional systems in \cite{Rybakowski82}.

Consider the space $\mathcal{X}$ of all topological spaces and the equivalence relation given by $Y\sim Z$ for $Y,Z\in \mathcal{X}$ if, and only if $Y$ is homotopy equivalent to $Z$, that is, there are continuous maps $f:Y\rightarrow Z$ and $g:Z\rightarrow Y$ such that $f\circ g$ and $g\circ f$ are homotopic to $id_Z$ and $id_Y$, respectively. Then, the quotient space $\mathcal{X}/ \sim$ describes the homotopy equivalent classes $[Y]$ of all topological spaces which have the same homotopy type. Intuitively, $[Y]$ describes all topological spaces which can be continuously deformed into $Y$.

Suppose $\Sigma$ is an invariant isolated set, that is, it is invariant with respect to positive and negative time of the semiflow, and it has a closed neighborhood $N$ such that $\Sigma$ is contained in the interior of $N$ with $\Sigma$ being the maximal invariant subset of $N$. 

Denote $\partial_e N\subset \partial N$ the \emph{exit set of $N$}, that is, the set of points which are not strict ingressing in $N$,
\begin{equation*}
    \partial_e N:=\{u_0\in N \text{ $|$ } u(t)\not\in N \text{ for all sufficiently small $t>0$}\} .
\end{equation*}

The \emph{Conley index} is defined as
\begin{equation*}
    C(\Sigma):=[N/ \partial_e N]
\end{equation*}
namely the homotopy equivalent class of the quotient space of the isolating neighborhood $N$ relative to its exit set $\partial_e N$. Such index is homotopy invariant and does not depend on the particular choice of isolating neighborhood $N$.

We compute the Conley index for two examples.

Firstly, the Conley index of a hyperbolic equilibria $u_+$ with Morse index $n$. Consider a closed ball $N\subset X$ centered at $u_+$ without any other equilibria in $N$, as isolating neighborhood. The flow provides a homotopy that contracts along the stable directions to the equilibria $u_+$. 
Then, $N$ is homotoped to a $n$-dimensional ball $B^n$ in the finite dimensional space spanned by the first $n$ eigenfunctions, related to the unstable directions. Note the exit set $\partial_e B^n=\partial B^n=\mathbb{S}^{n-1}$, since after the homotopy there is no more stable direction and the equilibria is hyperbolic. Therefore, the quotient of a $n$-ball and its boundary is an $n$-sphere,
\begin{equation*}
    C(u_+)=[N/\partial_e N]=[B^n/\partial_e B^n]=[B^n/\mathbb{S}^{n-1}]=[\mathbb{S}^n].
\end{equation*} 

Secondly, the Conley index of the union of two disjoint invariant sets, for example $u_-$ and $u_+$ with respective disjoint isolating neighborhoods $N_-$ and $N_+$. Then, $N_-\cup N_+$ is an isolating neighborhood of $\{ u_-,u_+\}$. By definition of the wedge sum 
\begin{align*}
    C(\{ u_-,u_+\})&=\left[\frac{N_-\cup N_+}{\partial_e (N_-\cup N_+)}\right]\\
    &=\left[\frac{N_-}{\partial_e N_-}\vee\frac{N_+}{\partial_e N_+}\right]=C(u_-)\vee C(u_+).
\end{align*}

The Conley index can be applied to detect heteroclinics as follows. Construct a closed neighborhood $N$ such that its maximal invariant subspace is the closure of the set of heteroclinics between $u_\pm$,
\begin{equation*}
    \Sigma=\overline{W^u(u_-)\cap W^s(u_+)}.
\end{equation*} 

Suppose, towards a contradiction, that there are no heteroclinics connecting $u_-$ and $u_+$, that is, $\Sigma=\{ u_-,u_+\}$. Then, the index is given by the wedge sum $C(\Sigma)=[\mathbb{S}^n]\vee [\mathbb{S}^m]$, where $n,m$ are the respective Morse index of $u_-$ and $u_+$.

If, on the other hand, one can prove that $C(\Sigma)=[0]$, where $[0]$ means that the index is given by the homotopy equivalent class of a point. This would yield a contradiction and there should be a connection between $u_-$ and $u_+$. Moreover, the Morse-Smale structure excludes connection from $u_+$ to $u_-$, and hence there is a connection from $u_-$ to $u_+$.

Hence, there are three ingredients missing in the proof: the Conley index can be applied at all, the construction of a isolating neighborhood $N$ of $\Sigma$ and the proof that $C(\Sigma)=[0]$.

\begin{pf}\textbf{of Proposition \ref{estabhetsquasi}}

($\impliedby$) This is also called liberalism in \cite{FiedlerRocha96}. Consider hyperbolic equilibria $u_-,u_+$ such that $i(u_-)=i(u_+)+1$ and satisfies both the Morse and the zero number permit conditions. Without loss of generality, assume $u_-(0)>u_+(0)$. 

It is used the Conley index to detect orbits between $u_-$ and $u_+$. Note that the semiflow generated by the equation \eqref{PDEquasi} on the Banach space X is admissible for the Conley index theory in the sense of \cite{Rybakowski82}, due to a compactness property that is satisfied by the parabolic equation \eqref{PDEquasi}, namely that trajectories are precompact in phase space. See Theorem 3.3.6 in \cite{Henry81}.

As mentioned above, in order to apply the Conley index concepts we need to construct appropriate neighborhoods and show that the Conley index is $[0]$. 

Consider the closed set
\begin{equation*}
    K({u_\pm}):= \left\{ u \in X \mid 
        \begin{array}{c} 
        z(u-u_-)=i(u_+)=z(u-u_+) \\
        u_+(0)\leq u(0)\leq u_-(0) 
        \end{array} \right\}
\end{equation*}

Consider also closed $\epsilon$-balls $B_\epsilon(u_\pm)$ centered at $u_\pm$ such that they do not have any other equilibria besides $u_\pm$, respectively, for some $\epsilon>0$. 

Define
\begin{equation*}
    N_\epsilon (u_\pm):=B_\epsilon(u_-)\cup B_\epsilon(u_+) \cup K({u_\pm}).    
\end{equation*}

The zero number blocking condition implies there are no equilibria in $K({u_\pm})$ besides possibly $u_-$ and $u_+$. Hence, $N_\epsilon (u_\pm)$ also has no equilibria besides $u_-$ and $u_+$. 

Denote $\Sigma$ the maximal invariant subset of $N_\epsilon$. We claim that $\Sigma$ is the set of the heteroclinics from $u_-$ to $u_+$ given by $\overline{W^u(u_-)\cap W^s(u_-)}$. 

On one hand, since $\Sigma$ is globally invariant, then it is contained in the attractor $\mathcal{A}$, which consists of equilibria and heteroclinics. Since there are no other equilibria in $N_\epsilon (u_\pm)$ besides $u_\pm$, then the only heteroclinics that can occur are between them.

On the other hand, Theorem \ref{Znuminvmfldquasi} implies that along a heteroclinic $u(t)\in\mathcal{H}$ the zero number satisfies $z^t(u-u_\pm)=i(u_+)$ for all time, since $i(u_-)=i(u_+)+1$. Therefore $u(t)\in K({u_\pm})$ and the closure of the orbit is contained in $N_\epsilon (u_\pm)$. Since the closure of the heteroclinic is invariant, it must be contained in $\Sigma$.

Lastly, it is proven that $C(\Sigma)=[0]$ in three steps, yielding the desired contradiction and the proof of the theorem. We modify the first and second step from \cite{FiedlerRocha96}, whereas the third remain the same. 

In the first step, a model is constructed displaying a saddle-node bifurcation with respect to a parameter $\mu$, for $n:=z(u_+-u_-)\in\mathbb{N}$ fixed,
\begin{equation}\label{prototype2}
    v_t=a(\xi,v,v_\xi)[v_{\xi\xi}+\lambda_nv]+ g_n(\mu,\xi,v,v_\xi)
\end{equation}
where $\xi\in [0,\pi]$ has Neumann boundary conditions, $\lambda_n=-n^2$ are the eigenvalues of the laplacian with $\cos(n\xi)$ as respective eigenfunctions, and 
\begin{equation*}
    g_n(\mu,\xi,v,v_\xi):=\left(v^2+\frac{1}{n^2}v^2_\xi-\mu\right)\cos(n\xi).
\end{equation*}

For $\mu>0$, a simple calculation shows that $v_\pm=\pm \sqrt{\mu}\cos(n\xi)$ are equilibria solutions of \eqref{prototype2} such that
\begin{equation}
    z(v_+-v_-)=n
\end{equation}
since the $n$ intersections of $v_-$ and $v_+$ will be at its $n$ zeroes. 

Moreover, those equilibria are hyperbolic for small $\mu>0$, such that $i(v_+)=n+1$ and $i(v_-)=n$. Indeed, parametrize the bifurcating branches by $\mu=s^2$ so that $v(s,\xi)=s\cos(n\xi)$, where $s>0$ correspond to $v_+$ and $s<0$ to $v_-$. Linearizing at the equilibrium $v(s,\xi)$ and noticing some terms cancel, the eigenvalue problem becomes
\begin{equation*}
    \eta v=a(\xi,s\cos(n\xi),sn\sin(n\xi))[v_{\xi\xi}+\lambda_n v]+\left[2s\cos(n\xi) v+\frac{2s\sin(n\xi) v_\xi}{-n}\right]\cos(n\xi)
\end{equation*}
where the unknown eigenfunction is $v$, corresponding to the eigenvalue $\eta$.

Hence $\eta_n(s)=2s$ is an eigenvalue with $v(s,\xi)$ its corresponding eigenfunction. 
Hence, by a perturbation argument in Sturm-Liouville theory, that is $\mu=0$ we have the usual laplacian with $n$ positive eigenvalues and one eigenvalue $\eta_n(0)=0$. Hence for small $\mu<0$, the number of positive eigenvalues persist, whereas for small $\mu>0$, the number of positive eigenvalues increases by 1. This yields the desired claim about hyperbolicity and the Morse index.

Now consider the quasilinear parabolic equation such that \eqref{prototype2} is the equilibria equation. The equilibria $v_\pm$ together with their connecting orbits form an isolated set
\begin{equation*}
    \Sigma_\mu(v_\pm):= \overline{W^u(v_-)\cap W^s(v_+)}
\end{equation*}
with isolating neighborhood $N_\epsilon(v_\pm)$, and the bifurcation parameter $\mu$ can also be seen as a homotopy parameter. Hence the Conley index is of a point by homotopy invariance as desired, that is,
\begin{equation}\label{conley1q}
    C(\Sigma_\mu(v_\pm))=C(\Sigma_0(v_\pm))=[0].
\end{equation}

In the second step, the $v_-$ and $v_+$ are transformed respectively into $u_-$ and $u_+$.

Recall $n=z(v_--v_+)=z(u_+-u_-)$. Hence, choose $\xi(x)$ a smooth diffeomorphism of $[0,\pi]$ that maps the zeros of $v_-(\xi)-v_+(\xi)$ to the zeros of $u_-(x)-u_+(x)$. Therefore, the zeros of $v_-(\xi(x))-v_+(\xi(x))$ and $u_-(x)-u_+(x)$ occur in the same points in the variable $x\in [0,\pi]$.

Consider the transformation 
\begin{align*}
    L: X &\to X   \\ 
    v(\xi)&\mapsto l(x)[v(\xi(x))-v_-(\xi(x))]+u_-(x)
\end{align*}
where $l(x)$ is defined pointwise through
\begin{align*}
    l(x):=
    \begin{cases}
        \frac{u_+(x)-u_-(x)}{v_+(\xi(x))-v_-(\xi(x))} &, \text{ if } v_+(\xi(x))\neq v_-(\xi(x)) \\
        \frac{\partial_x(u_{+}(x)-u_{-}(x))}{\partial_x(v_{+}(\xi(x))-v_{-}(\xi(x)))} &, \text{ if } v_+(\xi(x))= v_-(\xi(x))
    \end{cases}
\end{align*}
such that the coefficient $l(x)$ is smooth and nonzero due to the l'H\^opital rule. Hence, $L(v_-)=u_-$ and $L(v_+)=u_+$ as desired. Note we supposed $2\alpha+\beta>1$ so that solutions $u_\pm\in C^1$, hence $L$ is of this regularity as well. Moreover, $L$ is invertible with inverse having the same regularity. In particular, it is a homeomorphism, and hence a homotopy equivalence.

Moreover, the number of intersections of functions is invariant under the map $L$,
\begin{equation}
    z(L (v(\xi)-\tilde{v}(\xi)))=z(l(x)[v(\xi(x))-\tilde{v}(\xi(x))])=z(v(x)-\tilde{v}(x))
\end{equation}
and hence $K({v_\pm})$ is mapped to $K({u_\pm})$ under $L$. 

Consider $w(t,x):=L(v(t,\xi))$, hence the map $L$ modifies the equation \eqref{prototype2} into the following equation
\begin{equation}\label{IDK}
    w_t=\tilde{a}(x,w,w_x)w_{xx}+\tilde{f}(x,w,w_x)
\end{equation}
where the Neumann boundary conditions are preserved, and the terms $\tilde{a},\tilde{f}$ are
\begin{align*}
    \tilde{a}(x,w,w_x):=&\frac{x_\xi^2}{l(x)} \cdot a(x,L^{-1}(w),\partial_xL^{-1}(w))\\
    \tilde{f}(x,w,w_x):=&g_n(\mu,x,L^{-1}(w),\partial_xL^{-1}(w))+ (w_x\cdot x_{\xi\xi}-\partial_\xi^2u_-)-\frac{l_{\xi\xi} \cdot (w-u_-)_\xi}{l}\\
    &-\lambda a(x,L^{-1}(w),\partial_xL^{-1}(w))\cdot L^{-1}(w).
\end{align*}

Note that the equilibria $v_\pm$ are mapped into $w_\pm:=L(v_\pm)=u_\pm$, which are equlibria of \eqref{IDK}, with same zero numbers and Morse indices as $v_\pm$ and $u_\pm$.

The isolated invariant set $\Sigma_\mu(v_\pm)$ is transformed into $L(\Sigma_\mu(v_\pm))=\Sigma_\mu(w_\pm)$, which is still isolated and invariant, with invariant neighborhood $L(N_\epsilon(v_\pm))=N_\epsilon(w_\pm)$. Moreover, the Conley index is preserved, since $L$ is a homotopy equivalence,
\begin{equation}\label{conley2q}
   C(\Sigma_\mu(v_\pm))= C(L(\Sigma_\mu(v_\pm)))=C(\Sigma_\mu(w_\pm)).
\end{equation}

Hence, one identifies the equilibria $v_\pm$ in the model constructed \eqref{prototype2} with the equilibria $w_\pm=u_\pm$ from the equation \eqref{PDEquasi}, by preserving neighborhoods and the Conley index, since $L$ is a homotopy equivalence. The identified equilibria $u_\pm$ satisfy the equation \eqref{IDK}, and we still have to modify it to become \eqref{PDEquasi}. For such, we perform now a last homotopy between the solutions $w$ and $u$.

In the third step, we homotope the diffusion coefficient $\tilde{a}$ and nonlinearity $\tilde{f}$ from the equation \eqref{IDK} to be the desired diffusion $a$ and reaction $f$ from the equation \eqref{PDEquasi}. Indeed, consider the parabolic equation 
\begin{equation*}
    u_t=a^\tau(x,u,u_x)u_{xx}+f^\tau(x,u,u_x)
\end{equation*}
where
\begin{align*}
    a^\tau&:=\tau \tilde{a}+(1-\tau)a+\sum_{i=- \text{ , } +}\chi_{u_i}\mu_{u_i}(\tau)[u-{u_i}(x)]\\
    f^\tau(x,u,u_x)&:=\tau \tilde{f}+(1-\tau)f+\sum_{i=- \text{ , } +}\chi_{u_i}\mu_{u_i}(\tau)[u-{u_i}(x)]
\end{align*}
and $\chi_{u_i}$ are cut-offs beign 1 nearby $u_i$ and zero far away, the coefficients $\mu_i(\tau)$ are zero near $\tau=0$ and $1$ and shift the spectra of the linearization at $u_\pm$ such that uniform hyperbolicity of these equilibria is guaranteed during the homotopy. Note that $u_\pm$ have the same Morse indices, as solutions of both equations \eqref{PDEquasi} and \eqref{IDK}. Therefore, the $\mu_i(\tau)$ only makes sure none of these eigenvalues cross the imaginary axis. 

Consider $u_\pm$ and their connecting orbits during this homotopy, 
\begin{equation*}
    \Sigma^\tau(u_\pm):=\overline{W^u(u_-)\cap W^u(u_+)}.
\end{equation*}

Note that $\Sigma^\tau(u_\pm)\subseteq K({u_\pm})$, for all $\tau\in [0,1]$, since the dropping lemma holds throughout the homotopy. The equilibria $u_\pm$ do not bifurcate as $\tau$ changes, due to hyperbolicity. Choosing $\epsilon>0$ small enough, the neighborhoods $N_\epsilon(u_\pm)$ form an isolating neighborhood of $\Sigma^\tau(u_\pm)$ throughout the homotopy. Indeed, $\Sigma^\tau(u_\pm)$ can never touch the boundary of $K({u_\pm})$, except at the points $u_\pm$ by the dropping lemma. Once again the Conley index is preserved by homotopy invariance,
\begin{equation}\label{conley3q}
    C(\Sigma(u_\pm))=C(\Sigma^0(u_\pm))=C(\Sigma^\tau(u_\pm))=C(\Sigma^1(u_\pm))=C(\Sigma_\mu(w_\pm)).
\end{equation}

Finally, the equations \eqref{conley1q}, \eqref{conley2q} and \eqref{conley3q} yield that the Conley index of $\Sigma$ is the homotopy type of a point, and hence the desired result:
\begin{equation}\label{conley4q}
    C(\Sigma(u_\pm))=C(\Sigma_\mu(w_\pm))=C(\Sigma_\mu(v_\pm))=[0].
\end{equation}

    
\begin{flushright}
	$\blacksquare$
\end{flushright}
\end{pf}

\section{Example: Quasilinear Chafee-Infante}\label{sec:CIquasi}

In this section it is given an example of the theory above, namely, it is constructed the attractor of a quasilinear Chafee-Infante type problem, 
\begin{equation}\label{PDECIquasi}
    u_t= a(x,u,u_x)[u_{xx}+\lambda u[1-u^2]]
\end{equation}
where $n\in\mathbb{N}_0$, $x\in[0,\pi]$ has Neumann boundary conditions and initial data $u_0\in C^{2\alpha+\beta}([0,\pi])$ with $\alpha.\beta \in (0,1)$ and $a>0$, so that the equation generates a dynamical system in such space, as in \cite{Lunardi95}.

The equilibria equation describing the shooting curves is
\begin{align}\label{shootCIquasi}
\begin{cases}
    u_{\tau}&= p\\
	p_{\tau}&=-\lambda u[1-u^2]\\
	x_\tau&=1.
\end{cases}
\end{align}

Hence the shooting is exactly the same as the semilinear one. Thus, the permutation and the attractor are the same as the semilinear Chafee-Infante problem in \cite{FiedlerRocha96}. Therefore, both attractors are geometrically (connection-wise) the same. The only difference lies in the equilibria, and the parameter $\lambda$ must lie between two eigenvalues of the appropriate diffusion operator.

\section{Discussion}\label{sec:discuss}

We now mention two applications in the realm of general relativity and curve shortening flow, and lastly present two generalization proposals for the results in this chapter, namely a similar result for fully nonlinear parabolic equations, and fourth order parabolic equations.

It would be interesting to compute the attractor for other chosen nonlinearities $f$. In particular, one that has its application the construction of metrics at the event horizon of black holes, with a prescribed scalar curvature. For the case of self-similar solutions of the Schwarzschild metric type, the scalar curvature is chosen so that resulting parabolic equation is
\begin{equation*}
    u_t=u^2\left[u_{xx}+\frac{u_x}{\tan(x)}\right]-u(1-\lambda u^2)
\end{equation*}
where $\lambda\in\mathbb{R}_+$. Such problem was considered in \cite{FiedlerHellSmith}, where in such axially symmetric class it was shown that the equilibrium $u\equiv \lambda^{-1/2}$ bifurcates in an alternating sequence of pitchfork and transcritical bifurcations. Nevertheless, the attractor in such case was still not constructed due to the degenerate diffusion coefficients. Those attractors will be delt in the near future.



Another problem is that the diffusion coefficient $a(u)=u^2$ is not strictly parabolic when $u=0$. Either one has to construct Sturm attractors for parabolic equations with degenerate diffusion, or one restricts the phase space to $X^\alpha \cap \{ u>\epsilon \}$ for $\epsilon>0$ small, and hence obtain a subattractor within such subspace.

Moreover, the nonlinearity $f$ above does not satisfy the growth conditions that guarantees dissipativity. Numerical simulation of the shooting curve suggest that such nonlinearity is slowly nondissipative, as it is in \cite{RochaPimentel16}. We postpone this discussion, until all these tools have been sharpened and adapted for such case.

Another application is to obtain results in curve shortening flow. Considering a planar jordan curve flowing with respect to its curvature flow, then its curvature changes according to a quasilinear parabolic equation
\begin{equation*}
    k_t=k^2[k_{xx}+k]
\end{equation*}
as in \cite{Angenent91} and \cite{Angenent99}. Note there are several ways of modelling this same phenomena. 

Certain solutions are known to blow up in finite time. In particular, Grayson's theorem, in \cite{Grayson89}, guarantee that a strict convex curve shrinks to a point in finite time. The condition that a curve is strictly convex, mathematically is $k>\epsilon$ for $\epsilon>0$ fixed and guarantees strict parabolicity of the equation above. 

We believe that self-similar solutions of the type of the ODE blow up rate yields a dissipative nonlinearity, and the attractor of such equation is a single equilibria related to the circle. This would yield another proof of Grayson's theorem.

Lastly, we present two conjectures regarding attractors of parabolic equations.

The attractor construction in this chapter raises the question if one can construct the Sturm attractor for fully nonlinear second order parabolic and dissipative equation of the type
\begin{equation*}
    u_t=f(x,u,u_x,u_{xx})
\end{equation*} 
satisfying the parabolicity condition $\partial_q f>0$ for $q=u_{xx}$. This conjecture was already stated in \cite{Fiedler96}. Some results, such as shooting and obtaining permutations can be easily adapted. Others, like a Lyapunov functional and the full attractor construction are not so obvious. Those questions shall be treated in the near future. 

A last question is if it is possible to obtain any dynamical information of the attractor for higher order parabolic equations, such as the Swift-Hohenberg, 
\begin{equation*}
    u_t=-(1+\partial_x^2)^2u+\lambda u+f(u)
\end{equation*}
where $\lambda>1$. Such equation can be used to generate spatial-temporal patterns. For an overview on patterns and spirals, see \cite{Jia17}. 

This system has a Lyapunov functional. Hence, solutions either blow-up in finite time, or are global and then converge to equilibria. But the main tool from the second order equation in such case is not known, and probably does not hold: the dropping lemma. See \cite{PeletierTroy00}. 

The problem of rigorously constructing attractors for partial differential equations has started not so long ago, and still has a vast journey of discoveries and endless unanswered questions.

\medskip


\end{document}